\documentclass[11pt,a4paper]{article}
\usepackage{amsmath,amssymb,latexsym,graphics,epsfig}
\usepackage{hyperref}
\usepackage{amsthm}
\usepackage{graphicx,url}
\usepackage{amssymb}
\usepackage[utf8]{inputenc}
\usepackage[catalan]{babel}
\usepackage{color,soul}

\newcommand\blfootnote[1]{%
  \begingroup
  \renewcommand\thefootnote{}\footnote{#1}%
  \addtocounter{footnote}{-1}%
  \endgroup
}

\DeclareMathOperator{\sinc}{sinc}

\def\v{\mbox{\boldmath $v$}}

\def\vec0{\mbox{\boldmath $0$}}

\def\A{\mbox{\boldmath $A$}}

\def\D{\mbox{\boldmath $D$}}

\def\P{\mbox{\boldmath $P$}}

\def\T{\mbox{\boldmath $T$}}

\def\V{\mbox{\boldmath $V$}}

\def\Re{\mathbb R}

\begin{document}

\title{El mètode de les línies per a
la resolució num\`erica d'equacions en derivades parcials\\
The method of lines for numerical solutions of partial differential equations
\thanks{Aquesta recerca té el suport del projecte 2017SGR1087 de l'Agència de Gestió d'Ajuts Universitaris i de Recerca (AGAUR) del Govern de Catalunya. }}

\author{C. Dalfó\\
	\small{Departament  de Matemàtica, Universitat de Lleida, Igualada (Barcelona), Catalunya}\\
	\vspace{.5cm}
	\small{\url{cristina.dalfo@matematica.udl.cat}}\\
		M. A. Fiol\\
	\small{Deptartament de Matemàtiques, Barcelona Graduate School of Mathematics,}\\
		\small{Universitat Politècnica de Catalunya, Barcelona, Catalunya}\\
		\small{\url{miguel.angel.fiol@upc.edu}}
}

\date{}

\maketitle

\blfootnote{\begin{minipage}[l]{0.3\textwidth} \includegraphics[trim=10cm 6cm 10cm 5cm,clip,scale=0.15]{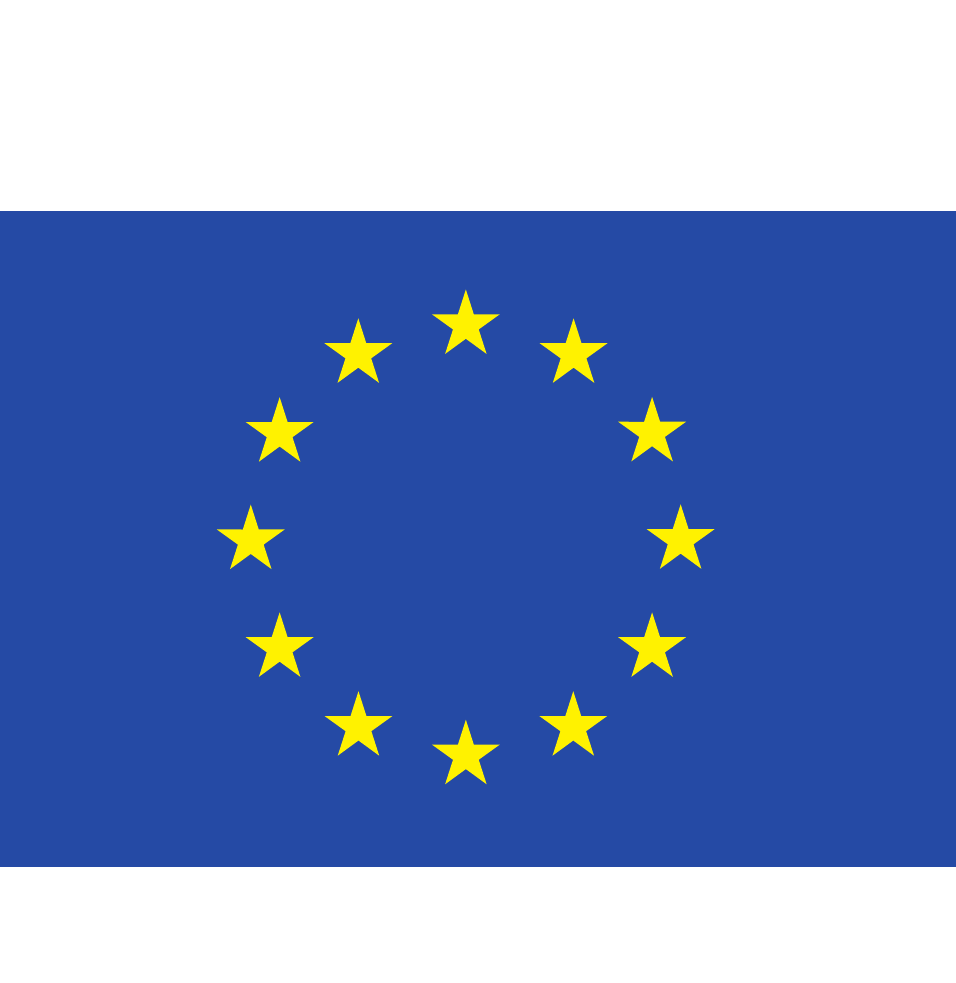} \end{minipage}  \hspace{-2cm} \begin{minipage}[l][1cm]{0.79\textwidth}
		The first author has received funding from the European Union's Horizon 2020 research and innovation programme under the Marie Sk\l{}odowska-Curie grant agreement No 734922.
\end{minipage}}

\begin{abstract}
En aquest treball es descriu un mètode numèric semi-discret per a la resolució d'un tipus d'equacions en derivades parcials. És conegut com el mètode de les línies (MOL, per les sigles en anglès),
i es basa en la discretització de totes les  variables involucrades, excepte una. I\l.lustrem l'aplicació del MOL amb la resolució de l'equació de Laplace en coordenades cartesianes. Els conceptes matemàtics que hi apareixen es comparen amb els corresponents al mètode analític de separació de variables. A més, mostrem que els resultats obtinguts amb el MOL  són aproximacions molt bones de les solucions analítiques.

\vskip.5cm

In this paper, we describe a semi-discrete method for a numerical resolution of a type of partial differential equations, called the method of lines (MOL). This method is based on the discretization of all but one of the variables of the problem. We illustrate this method by solving the Laplace equation in Cartesian coordinates. We compare the concepts used by the MOL with respect to the analytical method of variable separation. We show that the results obtained with the MOL are very good approximations of the analytical solutions.
\end{abstract}

\noindent{\em Paraules clau:} Equacions en derivades parcials, discretització d'una variable contínua, anàlisi numèric.

\section{Introducció}

El mètode de les línies, també anomenat mètode de semi-discretització per Zuazua \cite{Zu09a} és un mètode per a la resolució d'equacions en derivades parcials (EDP's) amb condicions de contorn i en regions amb una certa simetria. Es basa en discretitzar totes les variables, excepte una (si el problema ho permet, sovint es discretitzen les variables espacials i es manté contínua la variable temporal). Això condueix a un sistema d'equacions diferencials ordinàries (EDO's), les quals es poden resoldre pels mètodes numèrics habituals per a EDO's amb condicions inicials. Aquest mètode va ser introduït a principis dels anys 60 i, des de llavors, la seva exactitud i estabilitat han fet que es generessin una bona quantitat d'articles (vegeu, per exemple, Zafarullah \cite{Za70} i Verwer i Sanz-Serna \cite{vs84}).

Per tant, el MOL és un enfocament diferencial i de diferències finites per resoldre equacions en derivades parcials. Aquest mètode té moltes aplicacions en la resolució de problemes en Física i Enginyeria. Per exemple, aquest mètode s'ha usat àmpliament per experts en tècniques computacionals en problemes d'electromagnetisme, vegeu, per exemple, Berardia i Vurro \cite{BeVu16},  Diestel \cite{d84}, Sadiku i Obiozor \cite{SaOb}, o  Shakeri i Dehghan \cite{ShDe08}.

En aquest treball, il·lustrem el mètode de les línies amb l'equació de Laplace en el pla, per simplicitat. El procés a seguir consta dels punts següents:
\begin{enumerate}
\item
Discretització de l'equació diferencial en una direcció: Contracció de la regió de solucions a un conjunt de línies, i aproximació de les derivades corresponents.
\item
Obtenció d'un sistema d'equacions diferencials ordinàries acoblades.
\item
Transformació per obtenir un
sistema d'equacions diferencials ordinàries desacoblades.
\item
Resolució de les equacions diferencials.
\item
Transformació inversa i introducció de les condicions de contorn.
\item
Solució de les equacions.
\end{enumerate}

A més de les seves múltiples aplicacions, per exemple, com ja s'ha dit, en electromagnetisme, aquest mètode involucra moltes idees, conceptes i tècniques de la matemàtica discreta, juntament amb els mètodes bàsics de l'anàlisi diferencial. De fet, en el desenvolupament del mètode ens trobem amb els conceptes i els resultats següents:
les equacions diferencials en derivades parcials (e\l.líptiques, parabòliques i hiperbòliques), amb les seves aplicacions;
l'estudi dels 
problemes de contorn (amb condicions de Dirichlet o de Neumann);
les tècniques numèriques d'interpolació polinòmica i l'aproximació de derivades;
la teoria de
matrius tridiagonals semidefinides positives (les quals també apareixen en els grafs distància-regulars);
la teoria de
matrius circulants i la seva relació amb els grafs de Cayley sobre grups cíclics;
la diagonalizació de matrius i els algorismes del càlcul dels valors propis i vectors propis;
el teorema dels cercles de Gershgorin (per localitzar els valor propis);
els polinomis ortogonals de variable discreta (també importants en els grafs distància-regulars);
l'estudi de condicions d'ortogonalitat en el cas discret respecte del cas continu;
les transformacions discretes (canvis de base) com, per exemple, la Transformada Discreta de Fourier (DFT);
la resolució de les recurrències lineals i les diferències finites;
la resolució d'equacions diferencials ordinàries (EDO's) lineals;
l'aproximació i càlcul d'errors; i, finalment,
la representació gràfica i la interpretació de les solucions.

\section{El mètode de les línies (MOL): L'equació de Laplace}

	La transferència de calor a l'estat estacionari en sòlids segueix l'equació de Laplace i ha estat resolta directament per diversos mètodes numèrics, com el mètode successiu de sobre-relaxació, el mètode implícit de direcció alternativa i el mètode de transitoris falsos.
	Aquests mètodes són alternatius al mètode de les línies (MOL), el qual  
	és el que fem servir a continuació. El MOL consisteix en discretitzar totes les variables
	 excepte una. De vegades, el MOL també és anomenat mètode semi-analític (vegeu Subramanian i White \cite{sw04}).

\subsection{L'equació de Laplace en un domini rectangular}

\begin{figure}[t]
\begin{center}
\includegraphics[width=5cm]{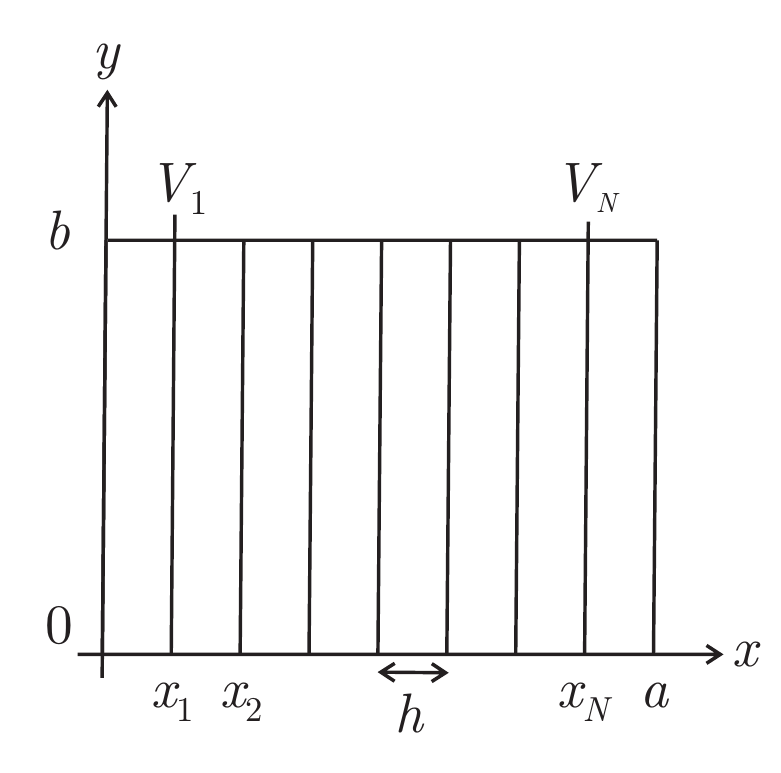}
\end{center}
\vskip-.5cm
\caption{Discretització de la variable $x$.}
\label{fig:dibuix1}
\end{figure}

Considerem l'equació de Laplace per a la funció potencial $V(x,y)$ en una regió rectangular $a\times b$.
\begin{align}
\label{Laplace1}
\nabla^2V &=\frac{\partial^2 V}{\partial x^2}+\frac{\partial^2 V}{\partial y^2}=0, \qquad x\in (0,a),\ y\in(0,b),
\end{align}
amb condicions de frontera que, per simplicitat, establim com
\begin{align}
V(0,y) &  =V(a,y)=0,\qquad\qquad\ \ y\in(0,b),\\
V(x,0) &  =0,\ V(x,b)=f(x), \qquad x\in(0,a).
\end{align}
El primer pas és discretitzar la variable $x\in(0,a)$ en els $N$ punts $x_j=j\Delta x$, per a $j=1,2,\ldots,N$, on el ``pas'' $h=\Delta x=\frac{a}{N+1}$ correspon a la distància entre dues de les $N$ línies resultants de la discretització, segons mostrem a la Figura~\ref{fig:dibuix1}.

Ara, iterem l'aproximació lineal de la derivada en dos punts successius,
$$
V'(x_j)\simeq \frac{V(x_{j+1})-V(x_{j})}{h},
$$
i obtenim l'anomenada ``aproximació centrada'' de la segona derivada
$$
V''(x_j)\simeq \frac{V'(x_{j})-V'(x_{j-1})}{h}\simeq
\frac{V(x_{j+1})-2V(x_{j})
+V(x_{j-1})}{h^2},
$$
amb error $|\epsilon_j|\le Ch^2$, per a una $C$ constant, que correspon a un mètode de segon ordre.
Amb notació simplificada,
\begin{equation}
\label{aproxV}
\frac{\partial^2 V}{\partial x^2}\simeq\frac{V_{j+1}-2V_{j}+V_{j-1}}{h^2},
\end{equation}
on $V_j$ denota la funció $V(x_j,y)$. Si substituïm \eqref{aproxV} en \eqref{Laplace1}, obtenim les $N$ equacions diferencials
\begin{equation}
\label{N edo's}
\frac{\partial^2 V}{\partial y^2}+\frac{1}{h^2}(V_{j+1}-2V_{j}+V_{j-1})=0,\quad ji=1,\ldots,N.
\end{equation}
Aleshores, si reemplacem $V(x,y)$ pel vector columna $\V=(V_1,\ldots,V_N)^\top$, podem escriure \eqref{N edo's} en forma matricial com
\begin{equation}
\label{sist-acoplado}
\frac{\partial^2 \V}{\partial y^2}-\frac{1}{h^2}\P\V=\vec0,
\end{equation}
on $\P$ és una matriu tridiagonal, és a dir,
\begin{equation}
\label{sist-acoplado-desarrollado}
 {\normalsize{\left(
     \begin{array}{c}
       \frac{\partial^2 V_1}{\partial y^2} \\
       \frac{\partial^2 V_2}{\partial y^2} \\
       \vdots \\
       \frac{\partial^2 V_{N-1}}{\partial y^2} \\
       \frac{\partial^2 V_N}{\partial y^2} \\
     \end{array}
   \right)-\frac{1}{h^2}
   \left(
     \begin{array}{ccccccc}
       p_l & -1 & 0 & \cdots & 0 & 0 & 0 \\
       -1 & 2 & -1 & \cdots & 0 & 0 & 0 \\
       0 & -1 & 2 & \cdots & 0 & 0 & 0 \\
        &  &  & \ddots &  & & \\
       0 & 0 & 0 & \cdots & -1 & 2 & -1 \\
       0 & 0 & 0 & \cdots & 0 & -1 & p_r \\
     \end{array}
   \right)
   \left(
     \begin{array}{c}
       V_1 \\
       V_2 \\
       \vdots \\
       V_{N-1} \\
       V_N \\
     \end{array}
   \right)=\vec0,
}}
\end{equation}
on els valors $p_l$ i $p_r$ depenen de les condicions de contorn a $x=0$ i $x=a$. Si la condició de contorn és de Dirichlet per l'esquerra [respectivament, per la dreta], aleshores $p_l=2$ [respectivament, $p_r=2$]; si la condició de contorn és de Neumann per l'esquerra [respectivament, per la dreta], aleshores $p_l=1$ [respectivament, $p_r=1$].

Per resoldre el sistema \eqref{sist-acoplado} d'equacions diferencials ordinàries (EDOs) lineals de segon ordre acoblades, necessitem transformar-lo en un sistema d'EDOs desacoblades.
Aleshores, la solució és diagonalitzar la matriu $\P$.
En aquest sentit, notem que $\P$ és simètrica i, per tant, els seus valors propis $\lambda_1,\ldots,\lambda_N$ són reals. A m\'es, el teorema de Gershgorin\footnote{Recordem que el Teorema de Gershgorin afirma que, donada una matriu $\A=(a_{jk})\in M_{n}(\mathbb{C})$, si es defineixen els cercles $D_{1},\ldots,D_{n}$ amb centre $a_{jj}$ i radi $r_{j}=\sum _{k\neq j}|a_{jk}|$, els valors propis de la matriu $A$ es troben en la unió dels $D_j$.} ens assegura la no negativitat, ja que han de complir $|\lambda_k-2|\le 2$, de manera que utilitzarem la notació $\lambda_k=\omega_k^2$ per a $k=1,\ldots,N$. De fet, com s'explica a la subsecció \ref{orto}, la matriu tridiagonal $\P$ té tots els valors propis diferents, i dona lloc a una successió de Sturm de polinomis ortogonals, els quals estan relacionats estretament amb els corresponents vectors propis (vegeu, per exemple, Godsil \cite[\S 8,5]{g93}).

Per tant, tenim que
$\T^\top \P\T=\D$, on $\T$ és la matriu ortogonal (amb columnes ortonormals corresponents als vectors propis de $\P$), és a dir,  $\T^{-1}=\T^\top$, i $\D$ és la matriu diagonal dels valors propis de $\P$. Aleshores, $\P=\T\D\T^{\top}$ que, si substituïm a \eqref{sist-acoplado}, ens dona
\begin{equation}
\label{sist-acoplado-transf1}
\nonumber
\frac{\partial^2 \V}{\partial y^2}-\frac{1}{h^2}\T\D\T^{\top}\V=\vec0,
\end{equation}
i, si multipliquem per $\T^{\top}(=\T^{-1})$,
\begin{equation}
\label{sist-acoplado-transf2}
\nonumber
\frac{\partial^2 \T^{\top}\V}{\partial y^2}-\frac{1}{h^2}\D\T^{\top}\V=\vec0,
\end{equation}
i, si utilitzem el {\em potencial transformat} $\overline{\V}:=\T^{\top}\V$, obtenim
\begin{equation}
\label{sist-acoplado-transf3}
\left(\frac{\partial^2}{\partial y^2}-\frac{1}{h^2}\D\right)\overline{\V}=\vec0.
\end{equation}
Per tant, obtenim $N$ equacions diferencials ordinàries (en la variable $y$), del tipus
$$
\overline{V}_j''=\frac{1}{h^2}\omega_j^2,\qquad j=1,\ldots,N,
$$
amb  solucions (el polinomi característic és
$s^2=\frac{\omega_j^2}{h^2}$,
amb arrels $s_{1,2}=\pm\frac{\omega_j}{h}$)
\begin{equation}
\label{soluciones}
 \overline{V}_j=\alpha_j e^{\frac{\omega_j}{h}y}+\beta_j e^{-\frac{\omega_j}{h}y}=A_j \cosh\frac{\omega_j y}{h}+B_j \sinh\frac{\omega_j y}{h},\quad j=1,\ldots,N.
\end{equation}
Per trobar el vector de potencials per a cada $j=1,\ldots,N$, només cal realitzar la {\em transformada inversa} de $\overline{\V}$, o sigui,
$$
\V=\T\overline{\V}.
$$
Finalment, imposem les condicions de contorn i resolem les equacions resultants per determinar les constants $\alpha_j,\beta_j$, o $A_j,B_j$.

\section{Diagonalització de la matriu $\P$ i ortogonalitat}

\subsection{Valors i vectors propis}

Abans de resoldre un exemple concret, en aquest apartat discutim la diagonalització de la matriu $\P$ i les propietats d'ortogonalitat involucrades.
Els resultats obtinguts, els valors propis $\lambda_k=\omega_k^2$ i les matrius de vectors propis $\T=(T_{jk})$, es mostren a la Taula \ref{tabla1}, segons les condicions de contorn a $x=0$ i $x=a$.

\begin{table}[t]
\begin{center}
\begin{tabular}{|c|c|c|c|}
  \hline
  Esquerra & Dreta & $T_{jk}$ & $\lambda_k=\omega_k^2$\\
  \hline\hline
  Dirichlet & Dirichlet & $\sqrt{\frac{2}{N+1}}\sin\left(\frac{jk\pi}{N+1}\right)$ & $4\sin^2\left(\frac{k\pi}{2N+2}\right)$\\
  \hline
  Dirichlet & Neumann & $\sqrt{\frac{2}{N+0.5}}\sin\left(\frac{j(k-0.5)\pi}{N+0.5}\right)$ & $4\sin^2\left(\frac{(k-0.5)\pi}{2N+1}\right)$\\
  \hline
  Neumann & Dirichlet & $\sqrt{\frac{2}{N+0.5}}\cos\left(\frac{(j-0.5)(k-0.5)\pi}{N+0.5}\right)$ & $4\sin^2\left(\frac{(k-0.5)\pi}{2N+1}\right)$\\
  \hline
  Neumann & Neumann & $\sqrt{\frac{2}{N}}\cos\left(\frac{(j-0.5)(k-1)\pi}{N}\right)$ si $k>1$ & $4\sin^2\left(\frac{(k-1)\pi}{2N}\right)$\\
          &         & $\frac{1}{\sqrt{N}}$ si $k=1$                                & \\
  \hline
\end{tabular}
\end{center}
\vskip-.25cm
\caption{Components de la matriu $T$ i valors propis $\lambda_k=\omega_k^2$ segons les condicions de contorn a $x=0$ i $x=a$.
}
\label{tabla1}
\end{table}

Per al càlcul procedim de la manera següent. L'equació a resoldre és
\begin{equation}
\label{ecuacion-autoval}
 \left(
     \begin{array}{ccccccc}
       p_l & -1 & 0 & \cdots & 0 & 0 & 0 \\
       -1 & 2 & -1 & \cdots & 0 & 0 & 0 \\
       0 & -1 & 2 & \cdots & 0 & 0 & 0 \\
        &  &  & \ddots &  & & \\
       0 & 0 & 0 & \cdots & -1 & 2 & -1 \\
       0 & 0 & 0 & \cdots & 0 & -1 & p_r \\
     \end{array}
   \right)
   \left(
     \begin{array}{c}
       v_1 \\
       v_2 \\
       \vdots \\
       v_{N-1} \\
       v_N
     \end{array}
   \right)= \lambda\left(
     \begin{array}{c}
       v_1 \\
       v_2 \\
       \vdots \\
       v_{N-1} \\
       v_N
     \end{array}
   \right),
\end{equation}
d'on
\begin{align}
(p_l-\lambda)v_1 &-v_2 =0, \label{eq1}\\
(2-\lambda)v_j &-v_{j-1} =v_{j+1},\quad j=2,\ldots,N-1 \label{eq2}\\
(p_r-\lambda)v_N &-v_{N-1} =0, \label{eq3}
\end{align}
on, per als valors usats de $p_l$ i $p_r$, el teorema de Gershgorin i el fet que la matriu és simètrica impliquen que els valors propis estan a l'interval $[0,4]$.

Per resoldre l'equació en diferències \eqref{eq2}, podríem usar la solució $v_j=s^j$, però s'obté un resultat més compacte amb $v_j=e^{ij\alpha}$ (és a dir, agafant $s=e^{i\alpha}=\cos \alpha+i\sin\alpha$). Si substituïm aquest valor a \eqref{eq2}, i dividim per $e^{i(j-1)\alpha}$ s'obté
\begin{equation}
e^{i\alpha}=(2-\lambda)-e^{-i\alpha}\quad \Rightarrow \quad 2-\lambda=2\cos \alpha.
\end{equation}
Per tant, l'equació característica és
\begin{equation}
\lambda=2(1-\cos\alpha)=4\sin^2\left(\frac{\alpha}{2}\right),
\end{equation}
cosa que significa que, per a cada valor propi $\lambda$, tenim dos possibles valors de $\alpha$:
$$
\alpha_{1,2}=\pm 2\arcsin\left(\frac{\sqrt{\lambda}}{2}\right).
$$
Així, amb $\alpha=\alpha_1$, la solució general de la recurrència és
\begin{equation}
\label{soluc-gen}
v_j=Ae^{ij\alpha}+Be^{-ij\alpha}.
\end{equation}
Per trobar els possibles valors de les constants $A$, $B$, i $\alpha$, imposem les ``condicions inicials'' \eqref{eq1} i \eqref{eq3}.
Suposem que tenim condicions de Dirichlet per la dreta i l'esquerra, és a dir, $p_l=p_r=2$. (Els altres casos es resolen de manera similar).
De la primera equació, s'obté que
\begin{align*}
(2-\lambda)v_1-v_2 &= 2\cos\alpha(Ae^{i\alpha}+Be^{-i\alpha})
-Ae^{i2\alpha}-Be^{-i2\alpha}=0\\
 & \Rightarrow\quad A+B=0,
\end{align*}
on hem aplicat $\cos\alpha=\frac{e^{i\alpha}+e^{-i\alpha}}{2}$ i $\sin\alpha=\frac{e^{i\alpha}-e^{-i\alpha}}{2i}$.
Per tant, $B=-A$ i, llevat d'una constant multiplicativa, podem suposar que
$$
v_j=\sin(j\alpha),\qquad j=1,\ldots,N.
$$
D'altra banda, de \eqref{eq3} s'obté
\begin{align*}
(2-\lambda)v_1-v_2 & =2\cos\alpha(Ae^{iN\alpha}+Be^{-iN\alpha})
-Ae^{i(N-1)\alpha}-Be^{-i(N-1)\alpha)}=0\\
 & \Rightarrow\quad \sin((N+1)\alpha)=0.
\end{align*}
Per tant, els possibles valors de $\alpha$ són $\alpha_k=\frac{k\pi}{N+1}$, per a $k=1,\ldots, N$, amb el que obtenim els valors propis
\begin{equation}
\label{autovalores}
\lambda_k=4\sin^2\left(\frac{k\pi}{2(N+1)}\right),\quad k=1,\ldots,N,
\end{equation}
els quals compleixen que $\lambda_1<\lambda_2<\cdots<\lambda_N$,
amb els respectius vectors propis
\begin{equation}
\label{autovectores}
\v_k=(v_{jk})=\left(\sin\left(\frac{k\pi}{N+1}\right), \sin\left(\frac{2k\pi}{N+1}\right),\ldots,\sin\left(\frac{Nk\pi}{N+1}\right)\right)^{\top}.
\end{equation}
Finalment, com que les columnes de $\T$ han de ser ortonormals, tenim
$$
\|\v_k\|^2=\sum_{j=1}^N \sin^2\left(\frac{jk\pi}{N+1}\right)=\frac{1}{2}\sum_{j=1}^N \left(1-\cos\left(\frac{2jk\pi}{N+1}\right)\right)=\cdots = \frac{N+1}{2},
$$
on hem usat que $\sin^2 x=\frac{1}{2}(1-\cos(2x))$ i la suma d'una sèrie geomètrica  $\sum_{j=1}^N r^j= \frac{r(1-r^N)}{1-r}$.
Concloem que els elements de la matriu $\T$ són
\begin{equation}
\label{T(jk)}
T_{jk}=(v_{jk})=\sqrt{\frac{2}{N+1}}\sin\left(\frac{jk\pi}{N+1}\right), \quad j,k=1,\ldots,N,
\end{equation}
com s'indica a la Taula \ref{tabla1}.

\subsection{Ortogonalitat}
\label{orto}

Considerem ara les propietats d'ortogonalitat que apareixen en el MOL, i les comparem amb les corresponents del mètode analític.
En primer lloc, notem que, per ser una matriu simètrica, els elements de les subdiagonals de la matriu $\P=(p_{jk})$ satisfan $p_{j,j+1}p_{j+1,j}>0$ per a $j=1,\ldots,N$ i, per tant,  donen lloc a una família de polinomis ortogonals $p_0,p_1,\ldots, p_{N-1}$ de variable discreta\footnote{És a dir, polinomis ortogonals respecte d'un producte escalar discret en els punts d'una malla $x_1,\ldots,x_N\in \Re$: $\langle f, g\rangle=\sum_{j=1}^N w_jf(x_j)g(x_j)$ per a uns certs ``pesos'' $w_j>0$.}. Aquest és també el cas, per exemple, de la matriu quocient d'un graf distància-regular, els polinomis ortogonals del qual corresponen als anomenats ``polinomis distància''. Els detalls d'aquesta aplicació, i altres de tipus combinatori, poden trobar-se a l'article de C\'{a}mara, F\`abrega, Fiol, i Garriga \cite{cffg09}. Aquests polinomis constitueixen una successió de Sturm (vegeu Godsil \cite[\S 8.5]{g93} or Chihara \cite{c78}) i satisfan una relació de recurrència de tres termes (com a l'expressió \eqref{eq2}) 
\begin{equation}
p_{k+1}(x)=(2-x)p_k(x)-p_{k-1}(x),\quad k=1,\ldots,N-2,
\label{eq-pol}
\end{equation}
inicialitzada amb $p_0(x)=1$ i $p_1(x)=2-x$. A més, $p_{N}(x)=(2-x)p_{N-1}(x)-p_{N-2}(x)$
correspon al polinomi característic de la matriu $\P$ i, per tant, les seves arrels són els valors propis $\lambda_k$, amb vectors propis
$$
\v_k=(p_0(\lambda_k),p_1(\lambda_k),\ldots,p_{N-1}(\lambda_k))^{\top},\qquad k=1,\ldots, N,
$$
trobats anteriorment a \eqref{autovalores} i \eqref{autovectores} (aquests últims llevat d'una constant multiplicativa).
L'ortogonalitat esmentada s'obté respecte d'un producte escalar discret en els punts de la malla
$\lambda_1,\lambda_2,\ldots, \lambda_N$, amb $\lambda_1<\lambda_2<\cdots <\lambda_N$. 
Per tant, un cop que normalitzem
els vectors $\v_k$, és a dir $\overline{\v}_k=\v_k/\|\v_k\|$, obtenim de nou els elements de la matriu $\T$ donats a \eqref{T(jk)}:
$$
T_{j,k}=(\overline{\v}_k)_j,\qquad j,k=1,\ldots,N,
$$
és a dir,
$$
T_{j,k}=\sqrt{\frac{2}{N+1}}\sin\left(\frac{jk\pi}{N+1}\right)
=\frac{p_{j-1}(\lambda_k)}
{\left(\sum_{j=0}^{N-1}p_j^2(\lambda_k)\right)^{\frac{1}{2}}},
\qquad j,k=1,\ldots,N.
$$
(Compareu amb les funcions pròpies que apareixen en el cas analític, vegeu la secció \ref{sec:exacte}).
En aquest context,
observem que el potencial transformat $\overline{\V}=\T^{\top}\V$ no és més que el vector de coeficients de Fourier obtingut a l'expressar $V$ en termes d'uns polinomis ortogonals (que, en els punts de la malla, tenen els valors donats per cada fila de $\T$).

\begin{figure}[t]
\begin{center}
\includegraphics[width=8cm]{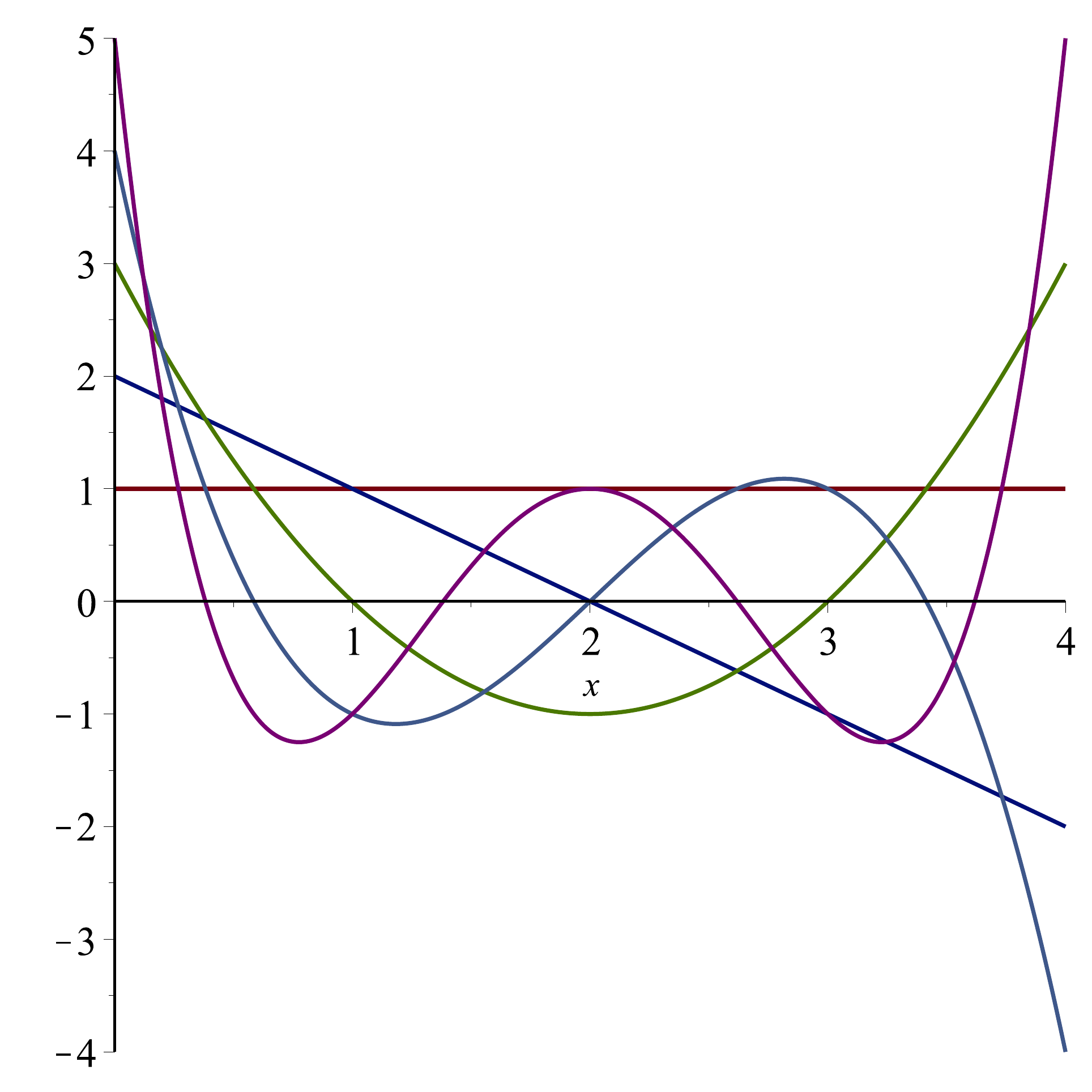}
\end{center}
\vskip-.5cm
\caption{Els polinomis ortogonals $p_j(x)$, per a $j=0,\ldots,4$.}
\label{fig:pols-ort}
\end{figure}

Per exemple, en el cas de la matriu $\P$ a \eqref{ecuacion-autoval}
amb $N=5$ i $p_l=p_r=2$, s'obtenen els polinomis
\begin{eqnarray*}
p_0(x)= &1,\\
p_1(x)= &2-x,\\
p_2(x)= &x^2-4x+3,\\
p_3(x)= &-x^3+6x^2-10x+4, \\
p_4(x)= &x^4-8x^3+21x^2-20x+5, \\
p_5(x)= &-x^5+10x^4-36x^3+56x^2-35x+6.
\end{eqnarray*}
Els valors propis de $\P$ són (arrels de $p_5(x)$):
$$
\lambda_1=2-\sqrt{3},\ \lambda_2=1,\ \lambda_3=2,\ \lambda_4=3,\ \lambda_5=2+\sqrt{3}.
$$
A la Figura \ref{fig:pols-ort} es mostren les gràfiques de $p_j(x)$, per a $j=0,\ldots,4$. Noteu les paritats (parells i senars) d'aquests polinomis respecte de $x=\lambda_3=2$, induïdes per les corresponents simetries de la matriu

$$
\T=\left(
\begin{array}{ccccc}
\sqrt{3}/6 & 1/2 & \sqrt{3}/3 & 1/2 & \sqrt{3}/6\\
1/2 & 1/2 & 0 & 1/2 & 1/2\\
\sqrt{3}/3 & 0 & -\sqrt{3}/3 & 0 & \sqrt{3}/3\\
1/2 & -1/2 & 0 & 1/2 & -1/2\\
\sqrt{3}/6 & -1/2 & \sqrt{3}/3 & -1/2 & \sqrt{3}/6\\
\end{array}
\right).
$$

\section{Un exemple}
\label{sec:exemple}

Ara anem a veure com apliquem el MOL a la resolució d'un problema concret.

\begin{figure}[t]
\begin{center}
\includegraphics[width=6cm]{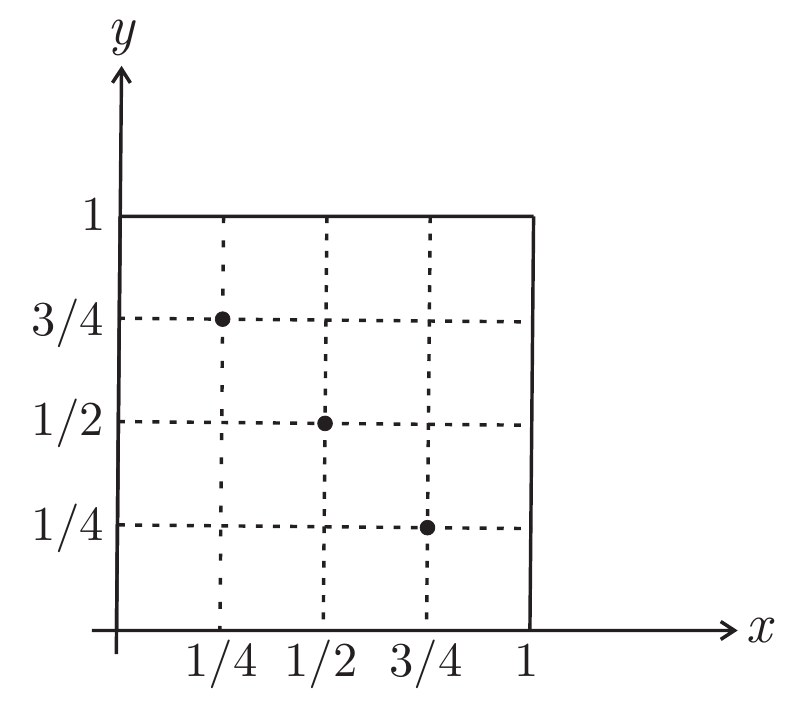}
\end{center}
\vskip-.5cm
\caption{Regió on volem resoldre l'equació de Laplace amb el MOL.}
\label{fig:dibuix2}
\end{figure}

Volem resoldre l'equació de Laplace (potencial)
$$
\frac{\partial^2 V}{\partial x^2}+\frac{\partial^2 V}{\partial y^2}=0,
$$
amb condicions de contorn de tipus Dirichlet $V(0,y)=V(a,y)=V(x,0)=0$ i $V(x,b)=100$, amb $a=b=1$ a la regió representada a la Figura~\ref{fig:dibuix2}. També volem calcular el valor del potencial en els punts indicats en aquesta figura.

La matriu $\P$ té valors propis (VAPS) $\omega_k^2=4\sin^2\frac{k\pi}{2N+2}$, per a $k=1,2,\ldots,N$, i la matriu de vectors propis (VEPS) $\T$ té components $T_{jk}=\sqrt{\frac{2}{N+1}}\sin\frac{jk\pi}{N+1}$ (vegeu la Taula~\ref{tabla1}).

Si escollim $N=15$, resulta que $h=\triangle x=\frac{1}{16}$, i les abscisses dels punts que volem calcular són $x_4=\frac{1}{4}$, $x_8=\frac{1}{2}$ i $x_{12}=\frac{3}{4}$.

De $V_j=\sum_{k=1}^N T_{jk} \overline{V}_k$, per a $i=1,2,\ldots,N$, obtenim la solució
$$
\left(
     \begin{array}{c}
       V_1 \\
       V_2 \\
       \vdots \\
       V_N\\
     \end{array}
   \right)=
   \left(
     \begin{array}{cccc}
       T_{11} & T_{12} & \cdots & T_{1N} \\
       T_{21} & T_{22} & \cdots & T_{2N} \\
          &  & \ddots &   \\
       T_{N1} & T_{N2} & \cdots & T_{NN} \\
     \end{array}
   \right)
   \left(
     \begin{array}{c}
       A_1 \cosh\frac{\omega_1 y}{h}+B_1 \sinh\frac{\omega_1 y}{h} \\
       A_2 \cosh\frac{\omega_2 y}{h}+B_2 \sinh\frac{\omega_2 y}{h} \\
       \vdots \\
       A_N \cosh\frac{\omega_N y}{h}+B_N \sinh\frac{\omega_N y}{h} \\
     \end{array}
   \right).
$$
Ara, per determinar les constants $A_j$ i $B_j$, imposem les condicions de contorn:
\begin{itemize}
\item
$V(x,0)=0 \Rightarrow$

$$
\left(
     \begin{array}{c}
       V_1 \\
       V_2 \\
       \vdots \\
       V_N\\
     \end{array}
   \right)=
   \left(
     \begin{array}{cccc}
       T_{11} & T_{12} & \cdots & T_{1N} \\
       T_{21} & T_{22} & \cdots & T_{2N} \\
        &  &   \ddots &  \\
       T_{N1} & T_{N2} & \cdots & T_{NN} \\
     \end{array}
   \right)
   \left(
     \begin{array}{c}
       A_1 \\
       A_2 \\
       \vdots \\
       A_N \\
     \end{array}
   \right)=\mathbf{0}
  $$
d'on
$A_j=0$, per a $j=1,2,\ldots,N$.
\item
$V(x,1)=100 \Rightarrow$

\begin{figure}[t]
	\begin{center}
		\includegraphics[width=10cm]{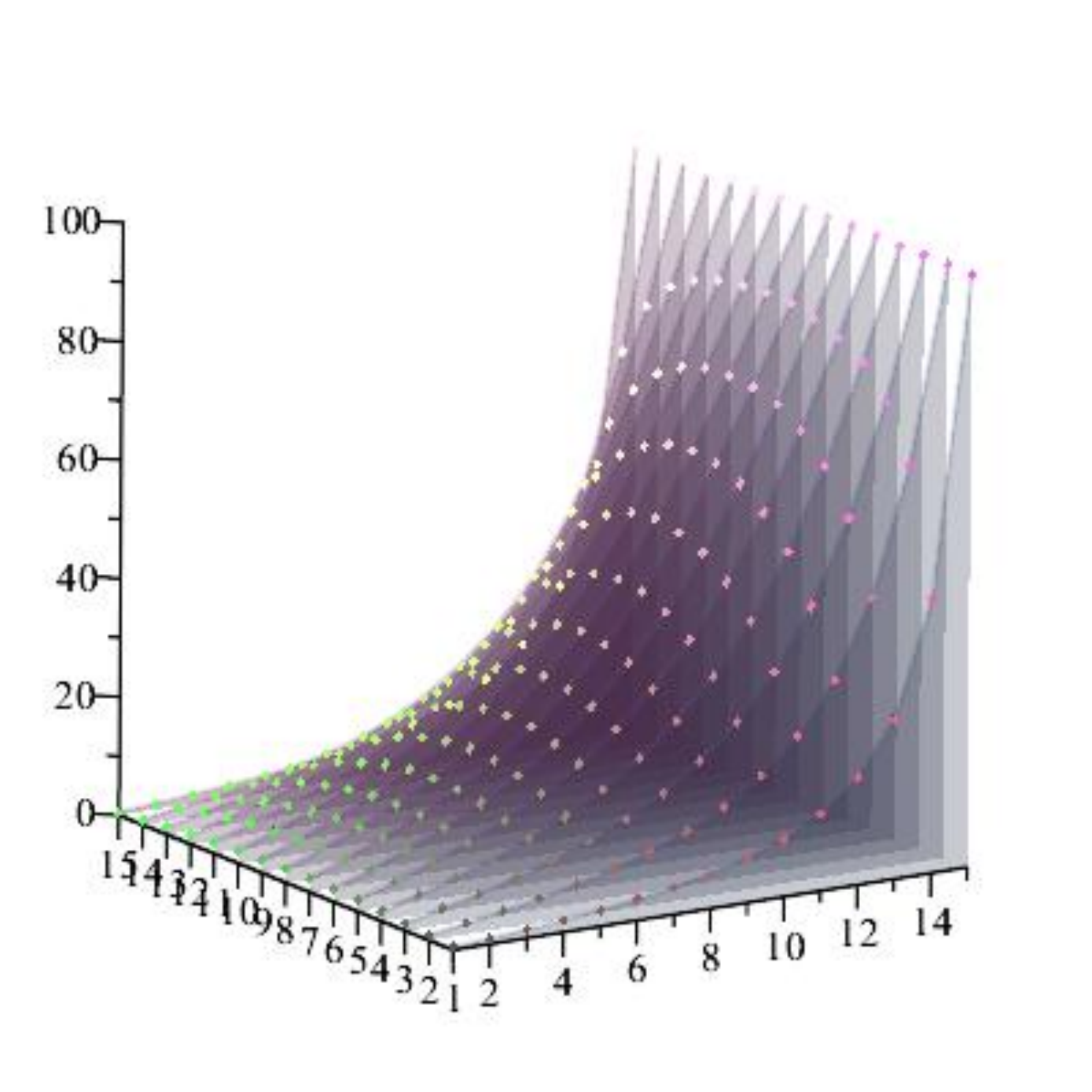}
	\end{center}
	\vskip-0.5cm
	\caption{Gràfica dels potencials $V_j(y)=V(x_j,y)$, amb $x_j=jh$ i $j=1,\ldots,15$, per a $V(x,1)=100$.}
	\label{fig:grafica-MOL}
\end{figure}

$$
\left(
     \begin{array}{c}
       V_1 \\
       V_2 \\
       \vdots \\
       V_N\\
     \end{array}
   \right)=
   \left(
     \begin{array}{cccc}
       T_{11} & T_{12} & \cdots & T_{1N} \\
       T_{21} & T_{22} & \cdots & T_{2N} \\
        &  &   \ddots &  \\
       T_{N1} & T_{N2} & \cdots & T_{NN} \\
     \end{array}
   \right)
   \left(
     \begin{array}{c}
       B_1 \sinh\frac{\omega_1}{h} \\
       B_2 \sinh\frac{\omega_2}{h} \\
       \vdots \\
       B_N \sinh\frac{\omega_N}{h} \\
     \end{array}
   \right)=
   \left(
     \begin{array}{c}
       100 \\
       100 \\
       \vdots \\
       100\\
     \end{array}
   \right),
  $$
d'on
  $$
  \left(
     \begin{array}{c}
       B_1 \sinh\frac{\omega_1}{h} \\
       B_2 \sinh\frac{\omega_2}{h} \\
       \vdots \\
       B_N \sinh\frac{\omega_N}{h} \\
     \end{array}
   \right)=
   \left(
     \begin{array}{cccc}
       T_{11} & T_{12} & \cdots & T_{1N} \\
       T_{21} & T_{22} & \cdots & T_{2N} \\
        &  &   \ddots &  \\
       T_{N1} & T_{N2} & \cdots & T_{NN} \\
     \end{array}
   \right)^\top
   \left(
     \begin{array}{c}
       100 \\
       100 \\
       \vdots \\
       100\\
     \end{array}
   \right):=
   \left(
     \begin{array}{c}
       c_1 \\
       c_2 \\
       \vdots \\
       c_N\\
     \end{array}
   \right).
  $$
Per tant,
  $$
  B_j=\frac{c_j}{\sinh\frac{\omega_j y}{h}},\qquad j=1,2,\ldots,N.
  $$
  \end{itemize}
Així, arribem a la solució
\begin{eqnarray*}
V_j(y) & = & V(x_j,y)=\sum_{k=1}^N T_{jk}B_k\sinh\left(\frac{\omega_j y}{h}\right)\\
       & = & \sum_{k=1}^N \sqrt{\frac{2}{N+1}}\sin\left(\frac{jk\pi}{N+1}\right)
       \frac{c_k}{\sinh\left(\frac{\omega_k}{h}\right)} \sinh\left(\frac{\omega_k y}{h}\right),
\end{eqnarray*}
la qual està representada a la Figura \ref{fig:grafica-MOL}.

En els punts demanats els valors que s'obtenen (calculats amb el programa de la Secció~\ref{sec:programes}) són
$$
\textstyle V(\frac{1}{4},\frac{3}{4})=43.1008,\quad V(\frac{1}{2},\frac{1}{2})=24.9644,\quad V(\frac{3}{4},\frac{1}{4})=6.7984
$$
que, comparats amb els valors calculats aplicant els resultats de la Secció~\ref{sec:exacte} també amb el programa de la Secció~\ref{sec:programes},
\begin{equation}
\label{valors-exactes}
\textstyle V(\frac{1}{4},\frac{3}{4})=43.2028,\quad  V(\frac{1}{2},\frac{1}{2})=24.9999,\quad V(\frac{3}{4},\frac{1}{4})=6.7971,
\end{equation}
demostren ser molt bones aproximacions.

\section{El mètode analític}
\label{sec:exacte}

Ara volem resoldre el problema anterior amb el mètode analític, per poder-lo comparar amb el resultat obtingut amb el MOL. 
És a dir, ja que podem conèixer la solució analítica, emprem l'equació de Laplace com a test pel mètode de les línies.  
	
\begin{figure}[t]
\begin{center}
\includegraphics[width=5cm]{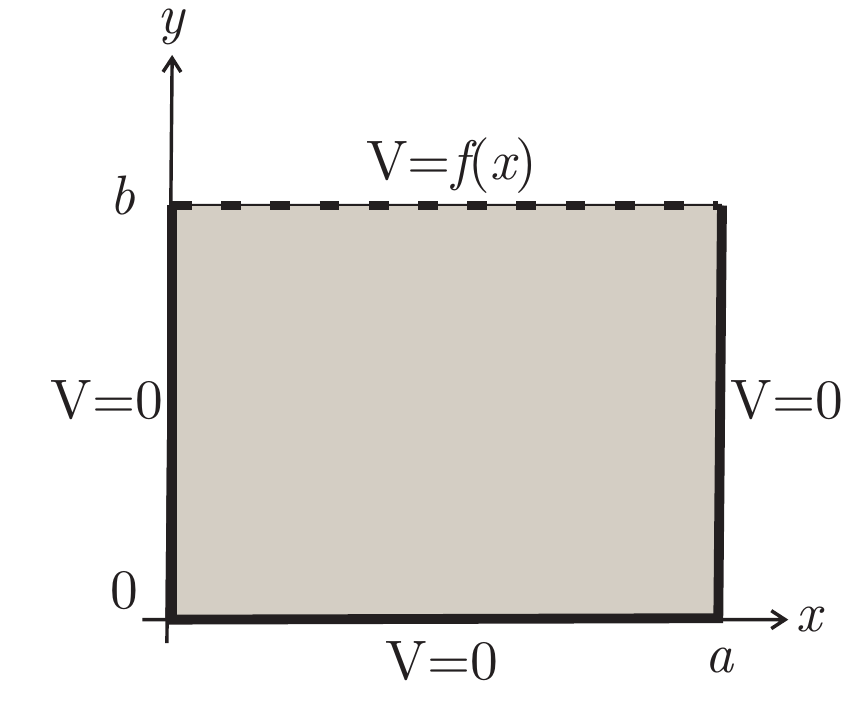}
\end{center}
\vskip-.5cm
\caption{Regió on volem resoldre l'equació de Laplace amb el mètode analític.}
\label{fig:dibuix3}
\end{figure}

Considerem doncs de nou l'equació
\begin{equation*}
\nabla^2V=\frac{\partial^2 V}{\partial x^2}+\frac{\partial^2 V}{\partial y^2}=0,
\end{equation*}
per a $x\in(0,a)$ i $y\in(0,b)$, amb condicions de contorn
\begin{align*}
V(0,y)&= V(a,y)=V(x,0)=0,\\
V(x,b)&=f(x),
\end{align*}
com es veu a la Figura \ref{fig:dibuix3}.

La solució analítica s'obté per separació de variables, és a dir, suposant que el potencial és de la forma
$V(x,y)=X(x)Y(y)$, d'on trobem que $X$ i $Y$ satisfan les equacions
\begin{align*}
X''(x)+\omega^2 X(x)  &=0,\\
Y''(y)-\omega^2 Y(y) &=0.
\end{align*}
Imposant les condicions de contorn homogènies, resulta que $\omega=\omega_k=\frac{k\pi}{a}$ per a $k=1,2,\ldots$ i obtenim  les funcions
\begin{align*}
X(x)  &=B\sin\left(\frac{k\pi}{a}x\right),\\
Y(y)  &=C\sinh\left(\frac{k\pi}{a}y\right),
\end{align*}
amb $B$ i $C$ constants.
Per tant, les solucions del problema homogeni són:
$$
V_k(x,y)=A_k\sinh\left(\frac{k\pi}{a}y\right)\sin\left(\frac{k\pi}{a}x\right).
$$
Considerem ara una solució en forma de sèrie
$$
V(x,y)=\sum_{k=1}^{\infty}A_k\sinh\left(\frac{k\pi}{a}y\right)
\sin\left(\frac{k\pi}{a}x\right).
$$

\begin{figure}[t]
	\begin{center}
		\includegraphics[width=10cm]{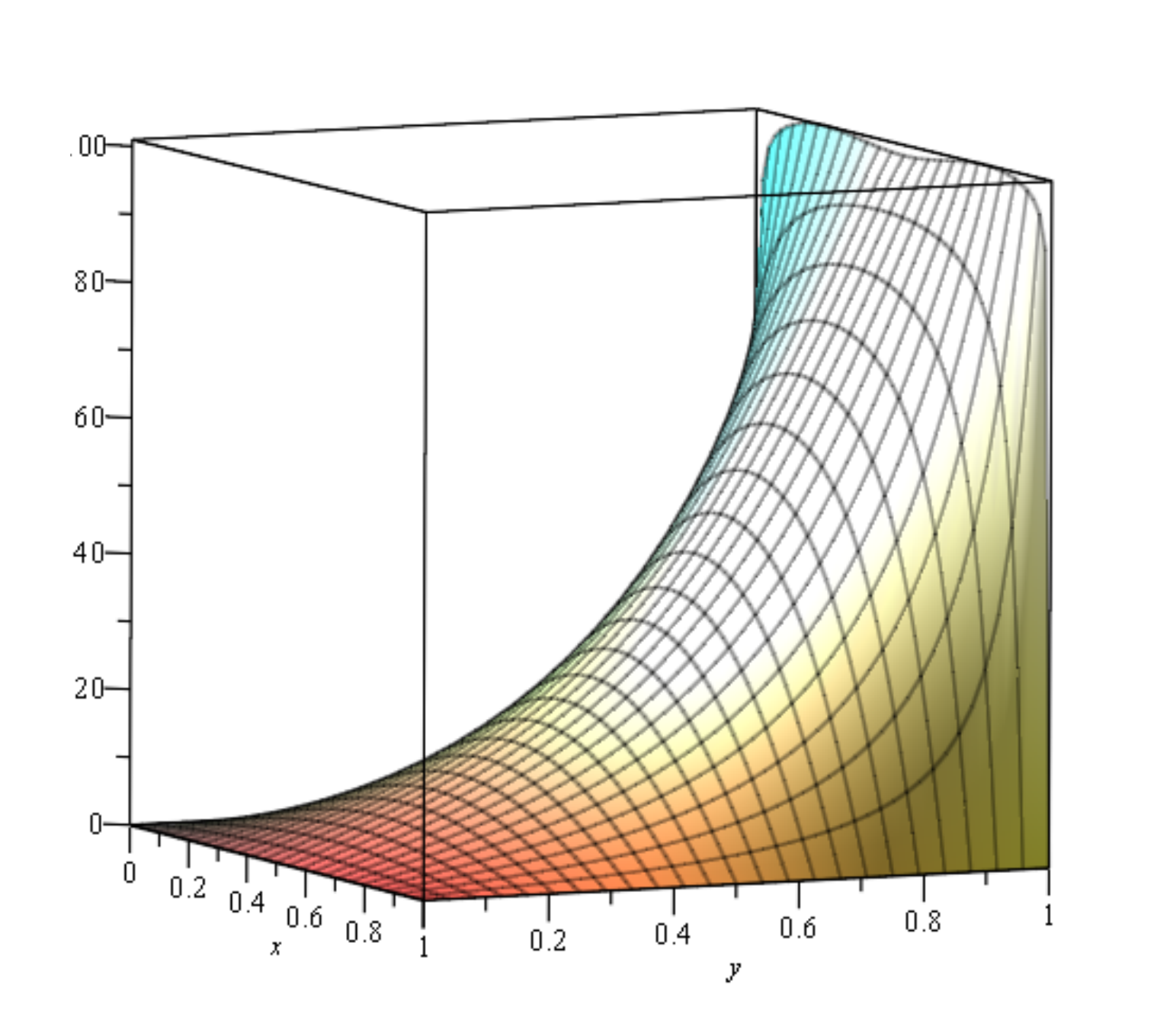}
	\end{center}
	\vskip-0.5cm
	\caption{Gràfica del potencial $V(x,y)$ en el cas $a=b=1$, $V(0,y)= V(a,y)=V(x,0)=0$ i $V(x,b)=100$.}
	\label{fig:grafica-exacte}
\end{figure}

Si imposem la condició de contorn no homogènia
$$
V(x,b)=f(x)=\sum_{k=1}^{\infty}A_k\sinh\left(\frac{k\pi}{a}b\right)
\sin\left(\frac{k\pi}{a}x\right),
$$
i considerem que $\{\sin\left(\frac{k\pi}{a}x\right):k=1,2,\ldots\}$ és una base ortogonal a $(0,a)$, obtenim el corresponent desenvolupament en sèrie de Fourier de $f(x)$, amb coeficients
$$
\alpha_k=A_k\sinh\left(\frac{k\pi}{a}b\right)=\frac{2}{a}\int_{0}^a f(x)\sin\left(\frac{k\pi}{a}x\right)\, dx.
$$
Per tant, la solució analítica és
$$
\displaystyle V(x,y)=\sum_{k=1}^{\infty}\alpha_k\frac{\sinh\left(\frac{k\pi}{a}y\right)}
{\sinh\left(\frac{k\pi}{a}b\right)}
\sin\left(\frac{k\pi}{a}x\right).
$$
Observem que les condicions de contorn \emph{iguals} a zero a $x=0$ i $x=a$ donen solucions periòdiques en \emph{sinus}, mentre que les condicions de contorn \emph{diferents} a $y=0$ i $y=b$ donen solucions no periòdiques en \emph{sinus hiperbòlics}.

En l'exemple anterior utilitzant el MOL, consideràvem el cas particular $a=b=1$ i $f(x)=100$. En aquest cas, ara obtenim que $\alpha_n=\frac{400}{n\pi}$ per a $n$ senar i  $\alpha_n=0$ per a $n$ parell. Aleshores, la solució analítica és
$$
\displaystyle V(x,y)=\frac{400}{\pi}\sum_{k=1}^{\infty}\frac{1}{2k-1}
\frac{\sinh((2k-1)\pi y)}
{\sinh((2k-1)\pi)}
\sin((2k-1)\pi x),
$$
la representació gràfica de la qual es mostra a la Figura \ref{fig:grafica-exacte}.

\section{Sobre l'error}
\label{sec:error}

Podem calcular l'error comès quan utilitzem el mètode de les línies, simplement amb la diferència entre els resultats obtinguts amb el MOL i el mètode analític. 
De fet, se sap que, sota condicions generals de regularitat de les solucions, l'error del mètode de les línies és de l'ordre de O$(h^2)$ (vegeu, per exemple, Schieser \cite{s91}, o Zuazua \cite{Zu09a,Zu09b}).

\begin{figure}[t]
\begin{center}
   \includegraphics[width=8cm]{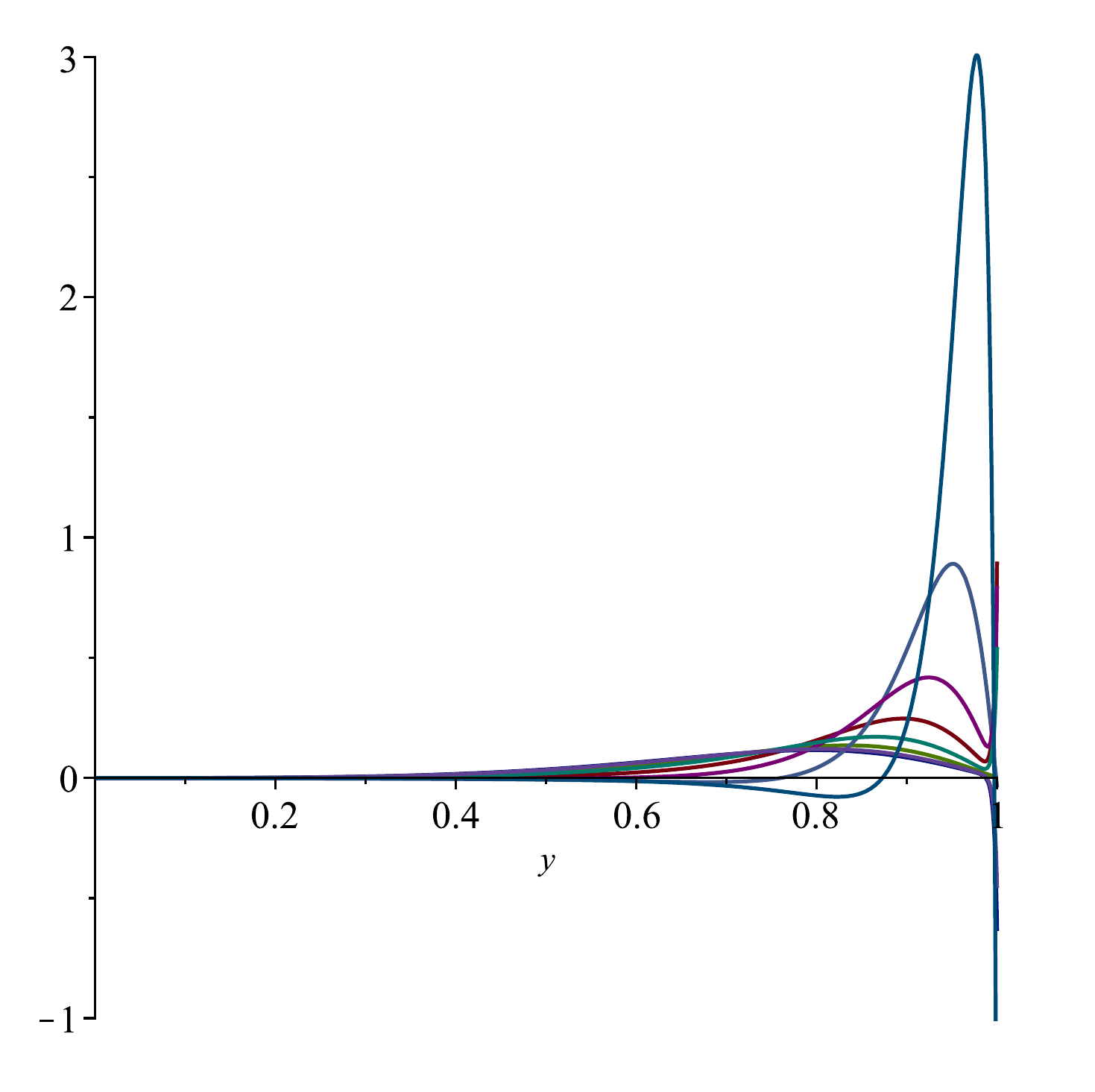}
\end{center}
   \vskip-0.5cm
   \caption{Gràfiques dels errors comesos en calcular el potencial amb el mètode de les línies en l'exemple de la Secció~\ref{sec:exemple} en comparació amb el mètode analític (Secció~\ref{sec:exacte}).}
\label{fig:grafica-error}
\end{figure}

Nosaltres estudiem primer el comportament de l'error
 en el cas de l'exemple anterior (Secció~\ref{sec:exemple}). Vegeu la Figura~\ref{fig:grafica-error}. Observem que tenim 8 funcions d'error (per a $x_j=j/16$, $j=1,\ldots,8$) malgrat que $N=15$ perquè, per simetria en les condicions de contorn, hi ha 7 funcions d'error que es superposen a d'altres 7 (i en queda una sense superposar). Com era d'esperar a causa de la discontinu\"{\i}tat en les condicions de contorn, els errors són més petits quan som més a prop del centre de la gràfica (és a dir, per a $x_8=1/2$).

\begin{figure}[t]
	\begin{center}
		\includegraphics[width=16cm]{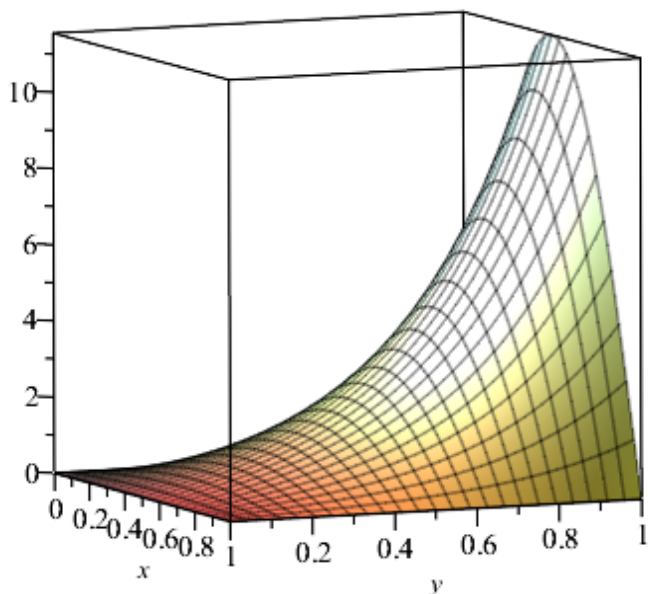}
	\end{center}
	\vskip-10cm
	\caption{Gràfica del potencial $V(x,y)$ en el cas $a=b=1$, $V(0,y)= V(a,y)=V(x,0)=0$ i $V(x,b)=\sinh(\pi)\sin(\pi x)$.}
	\label{fig:grafica-exacte-cas-especial}
\end{figure}

\begin{figure}[t]
	\begin{center}
		\includegraphics[width=16cm]{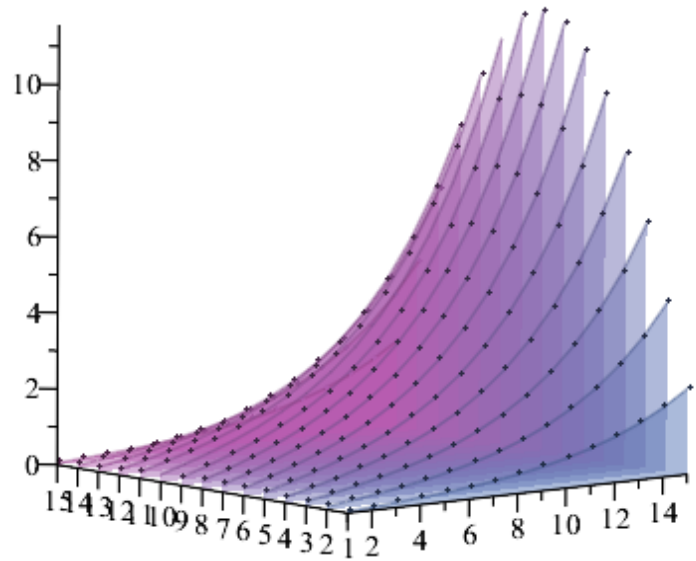}
	\end{center}
	\vskip-10cm
	\caption{Gràfica del potencial $V_j(y)=V(x_j,y)$, $x_j=jh$, $j=1,\ldots,15$, per a $V(x,b)=\sinh(\pi)\sin(\pi x)$.}
	\label{fig:grafica-aprox-cas-especial}
\end{figure}

Per obtenir resultats més precisos, 
considerem ara un altre exemple més senzill, el qual ens permet donar la solució numèrica amb una fórmula tancada. 
 Així, suposem ara que les condicions de contorn són
$V(0,y)=V(1,y)=V(x,0)=0$, com abans, i que $V(x,1)=\sinh(\pi)\sin(\pi x)$. Aleshores, per trobar la solució analítica no cal fer l'últim pas del desenvolupament en sèrie de Fourier, i s'obté 
\begin{equation}
\label{sol-exact-simple}
	V(x,y)=\sinh(\pi y) \sin (\pi x),
\end{equation}
	la gràfica de la qual es mostra a la Figura \ref{fig:grafica-exacte-cas-especial}.
Pel que fa al mètode de línies, comentem primer que aquest cas va ser resolt computacionalment per Subramanian i White \cite{sw04} utilitzant matrius exponencials, encara que sense arribar a una expressió tancada.
En canvi, les tècniques descrites aquí ens permeten obtenir la següent solució compacta
\begin{equation}
	\label{sol-aprox-simple}
	V(x_j,y)=\sinh[\pi\sinc(h) y]\sin(\pi x_j),
\end{equation}
on $x_j=jh$, $j=1,\ldots,N$, i $\sinc(x)=\frac{\sin(\pi x)}{\pi x}$. Vegeu la Figura \ref{fig:grafica-aprox-cas-especial}.
Notem que, com $\lim_{h\rightarrow 0}\sinc(h)=1$,  quan $h\rightarrow 0$, $j\rightarrow \infty$, i $jh\rightarrow x$, 
	i aleshores 
$V(x_j,y) \rightarrow V(x,y)$, com era d'esperar. De fet, en aquest cas és fàcil comprovar que, per a cada $x_j=x$ i $y$,  l'error $\epsilon=V(x,y)-V(x_j,y)$ és de l'ordre de $O(h^2)$, com ja s'ha comentat.

\section{Programa en MATLAB}
\label{sec:programes}

A continuació, escrivim en el llenguatge del programa MATLAB les instruccions per calcular el valor del potencial en els punts que es demanen a l'exemple anterior.

\newpage


\noindent $>>$ AA=1;\\
$>>$ BB=1;\\
$>>$ N=15;\\
$>>$ \%Determinem el vector alpha\\
$>>$ H=AA/(N+1);\\
$>>$ LAMBDA=2*sin((1:N)*pi*0.5/(N+1));\\
$>>$ ALPHA=LAMBDA/H;\\
$>>$ \%Calculem la matriu de transformació\\
$>>$ S=sqrt(2/(N+1));\\
$>>$ T=zeros(N,N);\\
$>>$ for J=(1:N)\\
\hspace*{1cm} for K=(1:N)\\
\hspace*{1.5cm}        T(J,K)=S*sin(J*K*pi/(N+1));\\
\hspace*{1cm}      end\\
\hspace*{.5cm}  end\\
$>>$ V=100*ones(N,1);\\
$>>$ C=inv(T)*V;\\
$>>$ A=ALPHA';\\
$>>$ B=C./sinh(BB*A);\\
$>>$ \%Calculem $V$ en els punts donats\\
$>>$ V1=0; V2=0; V3=0;\\
$>>$ for K=1:N\\
\hspace*{1cm} V1=V1+T(4,K)*B(K)*sinh(ALPHA(K)*0.75);\\
\hspace*{1cm} V2=V2+T(8,K)*B(K)*sinh(ALPHA(K)*0.5);\\
\hspace*{1cm} V3=V3+T(12,K)*B(K)*sinh(ALPHA(K)*0.25);\\
\hspace*{.6cm} end\\
$>>$ diary\\
$>>$ V1, V2, V3\\
\hspace*{1cm} V1 = 43.1008\\
\hspace*{1cm} V2 = 24.9644\\
\hspace*{1cm} V3 = 6.7984\\
$>>$ diary off

\vskip.5cm

%

\section{Altres exemples resolubles amb el MOL}

Finalment, detallem alguns altres exemples que es poden resoldre amb el mètode de les línies.

\begin{itemize}
\item[1.] 
L'equació de Laplace en la regió quadrada $1\times1$ amb condicions de contorn:
  \begin{itemize}
    \item[$(a)$] Dirichlet-Neumann.
    \item[$(b)$] Neumann-Dirichlet.
    \item[$(c)$] Neumann-Neumann.
  \end{itemize}
  (Els valors a utilitzar de la matriu $\T$ i els valors propis $\omega_k^2$ estan a la Taula \ref{tabla1}).
  \item[2.]
L'equació de Laplace en coordenades cilíndriques
  $\rho^2\frac{\partial^2 V}{\partial\rho^2}+\rho\frac{\partial V}{\partial\rho}+\frac{\partial^2 \phi}{\partial\phi^2}=0$, amb les  condicions de contorn $V(\rho,\phi)=0$ per a $0\leq\rho\leq a$ i
  $V(\rho,\phi)=V_0$ per a $a<\rho\leq b$, amb la discretització de la coordenada angular $\phi$.
La solució analítica és $V(\rho)=V_0\frac{\ln\left(\frac{\rho}{a}\right)}
{\ln\left(\frac{b}{a}\right)}$.
\item[3.]
L'equació d'ona $\frac{\partial^2 v}{\partial t^2}-\frac{\partial^2 v}{\partial x^2}=q(x,t)$, amb $0\leq x\leq L$ i $0\leq t\leq T$, les condicions inicials $v(x,0)=f_1(x)$ i $\frac{\partial v}{\partial t}(x,0)=f_2(x)$ per a $0\leq x\leq L$, la condició de contorn de Dirichlet $v(0,t)=g_1(x)$ i la condició no local $\int_0^L v(x,t)dx=g_2(x)$ per a $0\leq t\leq T$.
A la taula següent hi ha alguns casos particulars amb la corresponent solució analítica.
\begin{table}[t]
    \setlength{\tabcolsep}{1mm}
\small
\begin{center}
\begin{tabular}{|c|c|c|c|c|c|c|c|}
  \hline
   L & T & $q(x,t)$               & $f_1(x)$      & $f_2(x)$           & $g_1(t)$      & $g_2(t)$              & $v(x,t)$ analítica  \\
  \hline\hline
   1 & 4 & $0$                    & $0$           & $\pi\cos(\pi x)$   & $\sin(\pi t)$ & $0$                   & $\cos(\pi x) \sin(\pi t)$  \\
  \hline
   1 & 5 & $-2(x - t) e^{-x - t}$ & $0$           & $x e^{-x}$         & $0$           & $-2t e^{-t - 1}$      & $xt e^{-x-t}$ \\
         &   &                        & $0$           &                    & $0$           & $+ t e^{-t}$          & $xt e^{-x-t}$ \\
  \hline
   1 & 4 & $0$                    & $\cos(\pi x)$ & $0$                & $\cos(\pi t)$ & $0$                   & $\frac{1}{2}\cos(\pi (x+t))$   \\
            &   &                        &               &                    &               &                       & $+\frac{1}{2}\cos(\pi (x-t))$   \\
  \hline
   1 & 5 & $2x^5 + 2x^3 - 2x^2$   &            &  &           &  &  \\
            &   & $- (20x^3 + 6x - 2)\cdot$    &      $0$         &      $-x^5 - x^3 + x^2$              &   $0$            &      $\frac{t(t - 1)}{12}$                 & $(x^5 + x^3 - x^2)\cdot$   \\
            &   &  $(t^2 - t)$            &         &                    &               &                       & $(t^2 - t)$ \\
  \hline
\end{tabular}
\end{center}
\end{table}
\end{itemize}

\section{Conclusions}
	El mètode de les línies ha demostrat ser una bona eina per resoldre alguns tipus d'equacions en derivades parcials amb condicions de contorn. Com que es basa en discretitzar totes les variables excepte una, es pot utilitzar directament en l'equació de Laplace, la qual hem fet servir per il·lustrar aquest mètode. Una alternativa és aplicar el mètode de les falses transicions, afegint una derivada temporal i discretitzant les variables espacials.
	Si estem al pla $xy$, l'error d'aquest mètode alternatiu és de l'ordre de $\Delta_x^2+\Delta_y^2$. Per tant, és
	 millor aplicar el mètode semi-discret (o semi-analític), que és el que hem aplicat, i en el qual hem vist que l'error és de l'ordre de $\Delta_x^2$.
	
	Per una altra banda, 
	es pot 
	estudiar com canvia l'error si, en comptes de fer línies equiespaiades, fem línies situades als zeros dels polinomis de 
	Txebixov 
	(els quals tenen un millor comportament respecte de l'error). En aquest context, es pot consultar el treball de Youssef i Shukur \cite{yss14}. Finalment, proposem el següent problema obert: Estudiar  les semblances i les diferències entre el MOL i la integració al llarg de les rectes característiques de l'equació en derivades parcials de primer ordre $\frac{\partial V}{\partial t}+c\frac{\partial V}{\partial x}=f(x,t)$, amb $c>0$.
	
\section*{Agraïments}

Els autors volen agrair els comentaris i suggeriments del revisor anònim, els quals han contribuït considerablement a millorar aquest treball.


\end{document}